\documentclass{article}
\usepackage[utf8]{inputenc}

\usepackage{url}
\usepackage{optidef}
\usepackage{amsmath,amssymb,amsfonts,amsthm,enumerate}
\usepackage{eucal} 
\usepackage{xcolor}
\usepackage{youngtab}
\usepackage{young}
\usepackage{lscape}
\usepackage{tikz}
\usepackage{environ}
\usepackage{caption}
\usepackage{subcaption}
\usepackage{algorithm}
\usetikzlibrary{positioning, backgrounds}
\usepackage[noend]{algpseudocode}

\makeatletter
\newsavebox{\measure@tikzpicture}
\NewEnviron{scaletikzpicturetowidth}[1]{%
  \def\tikz@width{#1}%
  \begin{lrbox}{\measure@tikzpicture}%
  \BODY
  \end{lrbox}%
  \pgfmathparse{#1/\wd\measure@tikzpicture}%
  \BODY
}
\makeatother
\newcommand{\betti}{\beta}

\newcommand{\powset}{\mathcal{P}}

\newcommand{\witmap}{\Phi}
\newcommand{\cover}{\mathcal{U}}

\newcommand{\R}{\mathbb{R}}

\newcommand{\weightmap}{\pi}

\newcommand{\neighbors}{N}

\DeclareMathOperator{\Tr}{Tr}

\DeclareMathOperator{\ch}{ch}
\DeclareMathOperator{\conf}{Conf}

\DeclareMathOperator{\lazy}{Lazy}

\DeclareMathOperator{\skel}{Skel}

\DeclareMathOperator{\powcov}{PowCov}
\DeclareMathOperator{\cechplex}{\check{C}ech}
\DeclareMathOperator{\bary}{Sd}

\DeclareMathOperator{\nerve}{Nrv}

\DeclareMathOperator{\powdiag}{PowDiag}
\DeclareMathOperator{\alphaplex}{Alpha}

\DeclareMathOperator{\graph}{Graph}

\DeclareMathOperator{\delaunay}{Delaunay}

\DeclareMathOperator{\vor}{Vor}

\DeclareMathOperator{\im}{Im}
\DeclareMathOperator{\chull}{conv}

\DeclareMathOperator{\primalqp}{PrimalQP}
\DeclareMathOperator{\dualqp}{DualQP}

\newtheorem{thm}{Theorem}

\newtheorem{prop}{Proposition}
\theoremstyle{definition}
\newtheorem{defn}{Definition}

\newtheorem{example}{Example}

\makeatletter
\def\algbackskip{\hskip-\ALG@thistlm}
\makeatother

\usepackage{authblk}
\author[1]{Erik Carlsson}
\author[2]{John Carlsson}
\affil[1]{Department of Mathematics, UC Davis,
1 Shields ave. Davis, CA, 95618,
530-754-0274, 
ecarlsson@math.ucdavis.edu}
\affil[2]{Department of Industrial Engineering, University of Southern California, jcarlsso@usc.edu,
}

\date{}

\title{Computing the alpha complex using dual active set quadratic programming}

\begin{document}

\maketitle

\begin{abstract}

The alpha complex is a fundamental data structure from computational geometry, which encodes the topological type of a union of balls $B(x;r) \subset \mathbb{R}^m$ for $x\in S$,
including a weighted version that 
allows for varying radii.
It consists of the collection of 
``simplices'' $\sigma=\{x_0,...,x_k\} \subset S$, 
which correspond to nomempty $(k+1)$-fold intersections
of cells in a radius-restricted version of the
Voronoi diagram $\vor(S,r)$.
Existing algorithms for computing the alpha complex require that the points reside in low dimension because they begin by computing the entire Delaunay complex,
which rapidly becomes intractable,
even when the alpha complex is of a reasonable size.
This paper presents a method
for computing the alpha complex without computing 
the full Delaunay triangulation 
by applying Lagrangian duality, specifically an
algorithm based on dual quadratic programming 
that seeks to rule simplices out rather 
than ruling them in.

\end{abstract}

\section{Introduction}

Given a point cloud
and a threshold radius, the
alpha complex is a simplicial complex whose simplices correspond to
relationships between points that are relevant in the sense
that they are not too far apart. It is a generalization
of the Delaunay triangulation, another fundamental computational
geometric structure, which is the dual graph of the Voronoi
diagram. Alpha complexes are used in a diverse range of application
areas to study the shape of datasets, such as molecular biology \cite{liang1998analytical}, crystallography \cite{stukowski2014computational},
shape reconstruction \cite{1044617}, and persistent homology 
\cite{otter2015roadmap}.

\begin{figure}
\begin{centering}
\begin{subfigure}[b]{.32\textwidth}
\includegraphics[scale=.33]{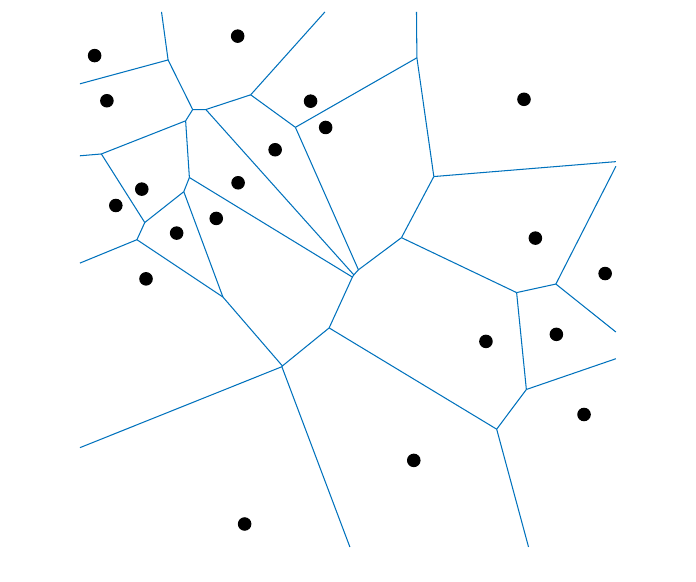}    
\caption{Voronoi diagram}
\label{fig:Voronoi-diagram}
\end{subfigure}
\begin{subfigure}[b]{.32\textwidth}
\includegraphics[scale=.33]{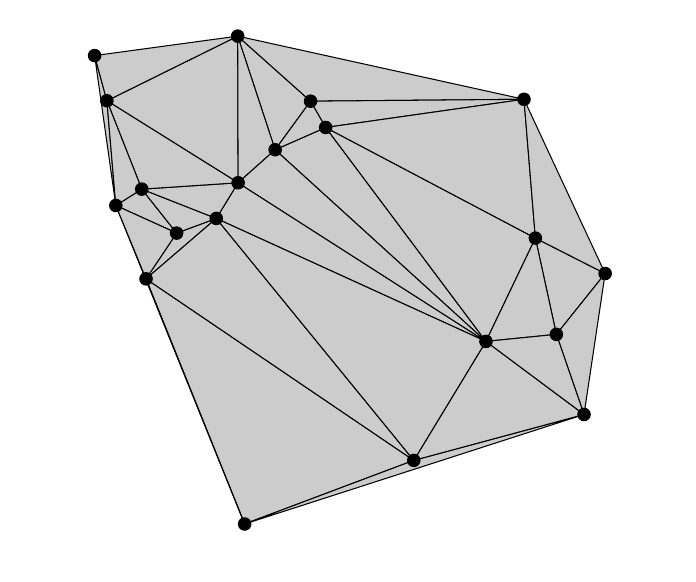}    
\caption{Delaunay triangulation}
\label{fig:Delaunay-triangulation}
\end{subfigure}
\begin{subfigure}[b]{.32\textwidth}
\includegraphics[scale=.33]{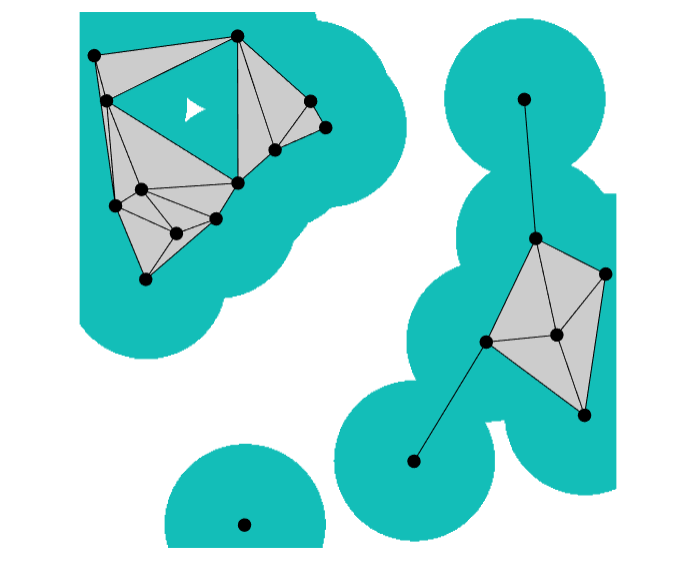}    
\caption{Alpha complex}
\label{fig:Alpha-complex}
\end{subfigure}
\end{centering}
\caption{The Voronoi partition associated
with a set of points (\ref{fig:Voronoi-diagram}), and the Delaunay
triangulation (\ref{fig:Delaunay-triangulation}). The alpha complex
(\ref{fig:Alpha-complex}) is a subset of the Delaunay triangulation
in which one removes simplices (which in the plane simply correspond
to points, edges, and triangles) that violate the distance property
that $V_{x}(r)\protect\neq\emptyset$.}
\label{fig:voronoi-delaunay-alpha}
\end{figure}

Formally, let $S=\{x_{1},...,x_{N}\}\subset\mathbb{R}^{m}$ be
a set of points, and let $r\geq 0$ be a nonnegative real number. The
\emph{radius-restricted Voronoi diagram} is the collection 
\[
\vor(S,r)=\{V_{x}(r):x\in S\}\,,
\]
where $V_{x}(r)=V_{x}\cap B(x;r)\subset\mathbb{R}^{m}$ is intersection
of the usual Voronoi cell 
\[
V_{x}=\{y\in\mathbb{R}^{m}:\|y-x\|\leq\|y-x'\|\,\forall x'\in S\}
\]
 with the ball $B(x;r)$ of radius $r$ centered about $x$. The alpha
complex is the \emph{nerve} of $\vor(S,r)$, that is, the simplicial
complex defined as the collection
\[
\alphaplex(S,r)=\left\{ \sigma\subset S:\bigcap_{x\in\sigma}V_{x}(r)\neq\emptyset\right\} .
\]
Simply put, a subset $\sigma=\{x_{i_{0}},\dots,x_{i_{k}}\}\subset S$
belongs to $\alphaplex(S,r)$ if there exists a point $y\in\mathbb{R}^{m}$
that is equidistant from every member of $\sigma$, i.e. $\rho:=\|y-x_{i_{0}}\|=\cdots=\|y-x_{i_{k}}\|\leq r$,
and furthermore, $\|y-x\|\geq\rho$ for all $x\in S$. A given subset
$\sigma$ of size $k+1$ is called a $k$-dimensional simplex 
of $\alphaplex(S,r)$.
The alpha complex is a subcomplex of $\delaunay(S)$, which is the
nerve of the full Voronoi diagram $\vor(S)$, and which agrees with
the Delaunay triangulation when the points of $S$ are planar points in general
position.  Figure \ref{fig:voronoi-delaunay-alpha}
shows one of these constructions. 

More generally, there
exists an extension of the alpha complex known as the weighted alpha
complex, in which the Voronoi cells are replaced by a power diagram,
which allows for balls of different radii. 
An example of a power diagram
and its associated weighted alpha complex is shown in Figure \ref{fig:nerve}
below. The weighted alpha complex gives rise to the \emph{alpha shapes}
\cite{edelsbrunner1992weighted},
whose applications include the aforementioned study of
biomolecules. A further extension is the wrap complex, which is used
in surface modeling \cite{edelsbrunner2003wrap,bauer2014morse}.

In terms of computational topology, the alpha complex is homotopy equivalent to the union of the cells 
\[
A=\bigcup_{x\in S}V_{x}(r)=
\bigcup_{x\in S}B(x;r)\subset\mathbb{R}^{m}.
\]
It can therefore be used to compute the topological type of 
a space which
can be covered by balls from the combinatorial data of which cells
intersect nontrivially. This is also true of the the \emph{\v{C}ech
complex} $\cechplex(S,r)$, defined as the nerve of the covering by
the balls as on the right, which satisfies 
\[
H_{*}(\alphaplex(S,r))\cong H_{*}(\cechplex(S,r)),
\]
both sides being isomorphic to $H_{*}(A)$. The alpha complex is by
definition a subcomplex of the \v{C}ech complex $\alphaplex(S,r)\subset\cechplex(S,r)$,
and it typically has far fewer simplices. This is advantageous, for
instance, for computing \emph{persistent homology} \cite{edelsbrunner2002topological, zomorodian2004computing},
noting
both the \v{C}ech and alpha complexes give rise to a family of complexes, which are ``filtered'' by varying the radius $r$.

Perhaps the most common construction for computing persistent homology is the Vietoris-Rips construction
\cite{hausmann1995vietoris}, especially through a highly efficient open source software tool
known as Ripser \cite{bauer2021ripser}.
\v{C}ech and alpha complexes can also be used for this purpose and have advantages over Vietoris-Rips in that they typically have far fewer simplices. Moreover,
because of theoretical guarantees stemming from the nerve theorem, their homology groups may be calculated exactly
from usual, non-persistent homology, 
which requires only Gaussian elimination.
One reason \v{C}ech and alpha complexes 
are not as commonly used is that Vietoris-Rips allows for non-Eucidean metrics, but a more crucial reason is the poor
scalability of the Delaunay construction in dimensions greater than three.

Most methods for computing the alpha complex begin by computing the
full Delaunay complex, and removing simplices which do not come from
a Voronoi face which has the restricted radius property \cite{cgal:dy-as3-23a,edelsbrunner1992weighted,edelsbrunner1994three,tralie2021cechmate}.
There are a wide range of highly efficient algorithms for computing
the Delaunay complex in dimensions $m\leq3$, many of which exploit
the empty circumsphere property, which states that a $3$-simplex
tetrahedron belongs to $\delaunay(S)$ in $\mathbb{R}^{3}$
if and only if its circumsphere contains no points \cite{bose2018flipping,hurtado1996flipping,aurenhammer1984optimal,klein1993randomized,bowyer1981computing,watson1981delaunay}.
In dimension $m>3$, one can still compute $\delaunay(S)$ by applying
flipping methods \cite{edelsbrunner1992incremental}, or by reducing
the problem to finding a convex hull in $\mathbb{R}^{m+1}$ as in
\cite{edelsbrunner1985voronoi}.
In terms of computing persistent homology, other authors have a combination of the alpha complex and the Vietoris-Rips construction to incorporate to improve efficiency \cite{mishra2023stability}.

In higher dimensions, it is often the case that the alpha complex
has a reasonable number of simplices, but the full Delaunay complex is far
too large to be computed, having on the order of $O(N^{\lceil m/2\rceil})$
simplices. In this situation, 
the general pipeline of the previous paragraph must
be replaced by one that does not compute the full complex. One algorithm that takes this into account is 
given in \cite{sheehy2015output}, whose complexity depends on
the total size of the output, and on bounds relating the pairwise
distance between points and the upper bound on the radius $r$.
 As an additional reduction, one is often only interested in subcomplex $\alphaplex_{\leq d}(S,r)$
consisting of simplices of dimension at most $d$.

A brute-force approach would be to formulate the existence of each
individual simplex $\sigma\in\alphaplex(S,r)$ as an optimization problem
\begin{mini}|l|
  {y\in \mathbb{R}^m}{\lVert
 y-x\rVert^2}{}{}
  \addConstraint{\lVert y-x_i\rVert}{\leq
 \lVert y-x_j\rVert}{\quad{}x_i \in \sigma,\ x_j \in S-\sigma}
   \addConstraint{\lVert y-x_i\rVert}{= \lVert y-x_j\rVert}{\quad{}x_i,x_j \in \sigma,\ i\neq j}
  \label{mini:qpdel}
 \end{mini} 
 where $x$ is any particular element of $\sigma$, all choices yielding the same
result. Specifically, we can conclude that a given simplex $\sigma$
is in $\alphaplex(S,r)$ when the constraints 
are feasible, and the minimizing value is at most $r^{2}$. 
By squaring the conditions, expanding and canceling terms, we actually see
that the above inequalities and equalities are in just linear constraints, making
\eqref{mini:qpdel} into a constrained (convex) quadratic program.

Computing
$\alphaplex_{\leq d}(S,r)$ would thus require solving
\[
\binom{N}{1}+\cdots+\binom{N}{d+1}
\]
such quadratic programs. In reality, many of these simplices
may be ruled out, including any simplex $\sigma$ which is not an
element of the \v{C}ech complex, or one whose faces have been determined
not to exist, assuming we are proceeding in order of increasing dimension.
Additionally, we only need to consider those constraints in (\ref{mini:qpdel})
coming from vertices $x_{j}$ which are neighbors in the one-skeleton
of $\cechplex(S,r)$. Despite these reductions, computing the alpha
complex directly is too burdensome to be practical for large values
of $N$ to be of practical value.

The point of this paper is to show that that this approach becomes
practical, provided that we use Lagrangian duality to solve (\ref{mini:qpdel})
instead of directly attacking the original primal problem. 
The fact that any
feasible point in the dual problem determines a lower bound for the
optimum of the primal problem is well-suited for this purpose
because it allows the
algorithm to terminate whenever a dual feasible point is found with
an objective value greater than $r^{2}$, which in many cases will
happen at an early stage.
Furthermore, the form of this particular dual problem, shown in
(\ref{maxi:dual}), has the property that the zero vector $\lambda=0$ is always feasible. Thus, there is no startup cost associated with identifying an initial feasible solution.
Another crucial benefit is that while the size  of the alpha complex is related to the rough dimension of the space traversed by the point
cloud $S$, there is essentially no dependence on the embedding, because dual programming algorithms depend only on the respective dot products.

Our method, which is straightforward to describe, is implemented for the more general weighted alpha complex, and is described in Algorithm \ref{alg:dualalpha}
below. Beyond using dual programming as described, we have taken advantage
of one further observation: the minimization problem (\ref{mini:qpdel})
for a given face is the same for that of the full cell $V_{x}$, except
that those inequalities determined by faces are replaced by the corresponding
equalities. The main loop of Algorithm \ref{alg:dualalpha} is written
in a way so that the the coefficients are only calculated once per
vertex, rather than once per potential simplex, which would otherwise
be a major computational cost. Algorithm \ref{alg:dualalpha} was
written in MAPLE, and is available at the first author's webpage:
\url{https://www.math.ucdavis.edu/~ecarlsson/}.
This includes an implementation of an elegant recent dual active set method
due to \cite{arnstrom2022dual}, which we used to solve the dual quadratic programs. 

In Section \ref{sec:examples}, we illustrate our algorithm in several examples which we validated using homology 
calculations, and which are also available online as MAPLE worksheets. We compared our answer against the output of persistence calculations which we carried out in Ripser. In some of those example there appears to be 
a potential computational advantage to using the alpha complex via our algorithm, as there often is for existing algorithms for the alpha complex in two or three dimensions \cite{somasundaram2021benchmarking}. 
However, the goal in comparing those answers is not to show a speed boost in persistent homology calculations, but rather to 
give a rigorous test of the correctness of the algorithm, which would fail to capture the correct homology if even a single simplex is incorrect. We make no comparison of the running time of our homology calculations, for which we used a general sparse matrix rank algorithm due to Dumas and Villard \cite{dumas2002sparse} instead of specialized methods. Intuitively, homology and persistent homology calculations of alpha complexes are expected to be faster than Vietoris-Rips calculations once the alpha complex has been computed, as the former is a subcomplex of the latter.

As we described above, the alpha complex has far reaching applications beyond persistent homology computations, firstly in that it produces concrete geometric models, which give rise to the alpha shapes. In terms of homology, it is also useful that the alpha complex provides exact answers rather that persistence diagrams, which we use in Section \ref{sec:config} to carry out an interesting calculation from geometric representation theory. Another recent application is due to the present authors, who discovered a hidden family of alpha complexes associated to the super-level sets of
an arbitrary kernel density estimator
in \cite{carlsson2023witness}. Implementing this construction in way that does not scale poorly with the embedding dimension was the motivation behind our main algorithm.

\subsection{Acknowledgments}

Both authors were supported by the
Office of Naval Research (ONR) N00014-20-S-B001 during this
project, which they gratefully acknowledge. 

\section{Preliminaries on computational topology}

We set some notation and background about filtered simplicial complexes,
refering to \cite{edelsbrunner2010computational} for more details.

\subsection{Computational topology}

\label{sec:comptop}

Let $S$ be a set of size $N$, which we assume is totally ordered.
\begin{defn}
\label{def:complex}
A simplicial complex $X$ on a vertex set
$S$ is a collection of nonempty subsets of $S$ which is closed under taking nonempty subsets.
\end{defn}
The elements of $X$ are called 
simplices, and are denoted $\sigma=[\sigma_{0},...,\sigma_{k}]$, 
using closed brackets to indicate that the elements are distinct and written in order. The number $k$ is called the dimension of $\sigma$, and the set of simplices of dimension $k$ is denoted $X_k$. 
If $\sigma\in X$ is a simplex, then the subsets
for which $\sigma'  \subset \sigma$ are called the faces of $\sigma$.
The collection of simplices of dimension at most $k$ 
is a subcomplex called the $k$-skeleton of $X$, which is denoted $\skel_k(X)$. 
For instance, $\skel_1(X)$ is a complex with only zero and one-dimensional simplices, which is the same data as a graph.

The following definition will also appear in Algorithm \ref{alg:dualalpha} below.
\begin{defn}
    \label{def:lazy}
Let $X$ be a complex. We define $\lazy_k(X)$ to
be the largest simplicial complex
$Z$ on the vertex set $S$ for which $\skel_k(Z)=\skel_k(X)$.
\end{defn} 
For instance, we have that $\lazy_{-1}(X)$ is the complete complex on the vertex set $S$, whereas $\lazy_0(X)$ is similar, but contains only those vertices which are in $X_0$, which need not be all of $S$. The one-dimensional lazy construction $\lazy_1(X)$ appears in 
the definition of the Vietoris-Rips and lazy Witness complexes \cite{desilva2008weak,hausmann1995vietoris,10.2312:SPBG:SPBG04:157-166},
which are widely used to compute persistent homology.

\begin{defn}
    \label{def:bary}
The Barycentric subdivision 
$\bary(X)$ is the complex whose vertices are the simplices in $X$, and whose $k$-simplices consist of strictly increasing flags 
$\sigma^{(0)}\subset \cdots \subset \sigma^{(k)}$
of elements of $X$.
\end{defn} 
\begin{defn}
    \label{def:geometric}
    If the vertex set $S$ is equipped with a map to $\mathbb{R}^m$, then the geometric realization is defined by
\begin{equation}
    \label{eq:geom}
    |X|=\bigcup_{\sigma \in X} \chull(\sigma) \subset \mathbb{R}^m
\end{equation}
where $\chull(\sigma)$ is the convex hull of the images of the vertices.
If no such map is given, the geometric realization is defined to be the standard one in which the $i$th element of $S$ is sent to the unit vector $e_i\in \mathbb{R}^N$.
\end{defn}
 Combining the two definitions, we see that a function $\witmap : X \rightarrow \mathbb{R}^m$ determines a linear map $|\bary(X)|\rightarrow \mathbb{R}^m$.

In persistent homology, one is interested in a nested family of complexes, depending on a 
real parameter $a$:
\begin{defn}
\label{def:filtcomplex}
A filtered complex is a pair
$(X,w)$ consisting of a simplicial 
complex $X$,
and a weight function 
$w:X\rightarrow \mathbb{R}$, which has the property that the subset $X(a)=w^{-1}(-\infty,a]\subset X$ is a complex for every $a$.
\end{defn}
The data of a pair $(X,w)$ and the corresponding nested collection of filtered complexes $X(a)$ and  are interchangeable. 
If $X$ is filtered by $w$, then there are induced filtrations on $X(a)$, $\skel_k(X)$, and $\lazy_k(X)$.
The filtration on $\lazy_k(X)$ is the one for which
$w(\sigma)$ is the max of $w(\sigma')$ as $\sigma'$ ranges over all elements of $\skel_k(X)$ which are faces of $\sigma$. The first two are simply by restriction.

If $X$ is a complex then the chain group is the 
vector space of all formal linear combinations
\begin{equation}
C_k=\left\{
\sum_{\sigma\in X_k} c_{\sigma} \sigma
\right\}
\end{equation}
where we will always take coefficients to be
elements of a finite field $c_{\sigma} \in \mathbb{F}_p$ for $p$ a prime. The $k$th homology group is given by
\begin{equation}
\label{eq:homology}
H_k(X,\mathbb{F}_p)=Z_k/B_k=\ker \partial_k/\im \partial_{k+1}
\end{equation}
where $\partial_k:C_{k}\rightarrow C_{k-1}$ is
defined on each basis vector $\sigma=[x_0,...,x_k]$ by
\[\partial_k(\sigma)=\sum_{i=0}^k (-1)^i 
[x_{0},...,\widehat{x_i},...,x_k]\]
The $k$th Betti number $\beta_k(X,\mathbb{F}_p)$ is 
the dimension of $H_k(X,\mathbb{F}_p)$.

Filtered complexes have the additional structure of persistent homology groups, which assemble the individual homology groups $H_k(X(a))$ for each value of $a$ into a family of 
filtered homology groups. Instead of individual Betti numbers, one has a collection of persistence intervals often called a barcode diagram, such as the Ripser outputs shown in Section \ref{sec:examples}. Roughly speaking, one can infer the Betti numbers of a point cloud by counting the significant intervals in that diagram. For an introduction to persistent homology, we refer to
\cite{edelsbrunner2010computational}.

\subsection{The Nerve theorem}

\label{sec:nerve}

Let $\mathcal{U}=\{U_x:x\in S\}$
be a collection of subsets of $\mathbb{R}^m$ 
with index set $S$.
\begin{defn}
 The \v{C}ech nerve, written
$\nerve(\mathcal{U})$, is the complex with vertex set $S$, and for which
\begin{equation}
    \label{eq:nerve}
    [x_0,...,x_k]\in X
    \Longleftrightarrow
    U_{x_0} \cap \cdots \cap U_{x_k} \neq \emptyset
\end{equation}
\end{defn}
Then the classical nerve theorem of Leray \cite{Leray1950LanneauSE} states:
\begin{thm}[Leray]
Suppose that $\mathcal{U}$ 
has the property that any $k$-fold intersection of the 
$U_x$ is contractible
(which occurs, for instance, if every $U_x$ is convex).
Then $\nerve(\mathcal{U})$ is homotopy 
equivalent to the union 
$\bigcup \cover \subset \mathbb{R}^m$. 
\end{thm}
In some cases, the nerve equivalence has an explicit form. If $\mathcal{U}$ is a covering by balls, then the linear map $|\nerve(\mathcal{U})|\rightarrow \bigcup \mathcal{U}$ determined by sending each vertex to the corresponding center induces the nerve equivalence. More generally, if every $U_x$ is convex, and we select any representatives $x_\sigma \in U_{\sigma_0} \cap \cdots \cap U_{\sigma_k}$, we have an induced map
$\Phi:|\bary(X)|\rightarrow \bigcup \mathcal{U}$, which also induces the nerve equivalence equivalence via the equivalence of $|X|$ with 
the subdivision $|\bary(X)|$. See \cite{bauer2023unified} for a proof, and more on the general setup of nerve theorems.

\subsection{Power diagrams}

\label{sec:powdiag}

Suppose that $S\subset \R^m$, and
let $p:S\rightarrow \mathbb{R}$ 
be a function, called the weight map. 
We now have a function
$\weightmap : \R^m \rightarrow \R$ given by
\begin{equation}
    \label{eq:weightmap}
\weightmap(y)=\min_{x\in S} \pi_x(y),\quad \pi_x(y)=\lVert y-x\rVert^2-p(x).    
\end{equation}
\begin{defn}
Let $S,p$ be as above and let $a\in \mathbb{R}$. Then the weighted ball cover denoted $\mathcal{U}=\powcov(S,p,a)$
is given by $\mathcal{U}=\{U_x:x\in S\}$,
where
\begin{equation}
    \label{eq:powcov}
U_x=\left\{y\in \mathbb{R}^m:
\pi_x(y)\leq a\right\}
\end{equation}
is either a closed ball, or is empty.
\end{defn}
\begin{defn}
The weighted power diagram is the covering $\powdiag(S,p,a) =\{U_x\cap V_x:x\in S\}$ where 
\begin{equation}
\label{eq:powdiag}V_x=\left\{y:\pi_x(y)\leq \pi_{x'}(y) \mbox{ for all $x'\in S$}\right\},
\end{equation}
and $U_x \in \powcov(S,p,a)$ ranges over the corresponding elements in the weighted ball cover.
\end{defn}

The \v{C}ech and alpha complexes are the filtered complexes which are the nerves of the weighted ball and power cover, resepectively:
\begin{defn}
\label{def:cech}    
The \v{C}ech complex denoted $(X,w)=\cechplex(S,p)$ is the
filtered complex determined by 
$X(a)=\nerve(\powcov(S,p,a))$. We let $\cechplex(S,p,a_1)$ 
be the filtered subcomplex $X(a_1)$ which is cut off at weight $a_1$.
\end{defn}
\begin{defn}
\label{def:alphaplex}
The weighted alpha complex
$(X,w)=\alphaplex(S,p)$ is the filtered complex for which
$X(a)=\nerve(\powdiag(S,p,a))$, with a similar definition of $\alphaplex(S,p,a_1)$.
\end{defn}
Said another way, we have a simplex 
$\sigma=[x_0,...,x_k] \in X(a)$
If the weighted Voronoi face $V_\sigma=V_{\sigma_0}\cap \cdots \cap V_{\sigma_k}$ is nonempty, and there exists a point $x\in V_\sigma$ satisfying $\pi_{x_i}(x)\leq a$ for any $i$, noticing that the $\pi_{x_i}$ all become equal when restricted to $V_\sigma$. Adopting the terminology of the witness complex, such a point $x$ is called a \emph{witness} for $\sigma$ because its existence determines that 
$\sigma \in X(a)$. Since $\pi_{x_i}$ is a quadratic function and $V_\sigma$ is convex, we have a unique minimizer 
$x_\sigma$ for every $\sigma \in X(a)$. The collection of these points is described as a map:
\begin{defn}
    \label{def:witmap}
Let $X=\alphaplex(S,p)$. The 
\emph{witness map} is the 
function $\Phi:X \rightarrow \mathbb{R}^m$ 
which carries each simplex
$\sigma=[x_{0},...,x_{k}]$ to the unique element
$x_\sigma \in V_{\sigma}$ that minimizes $\pi_{x_j}(x)$, 
which is independent of the choice of $j$.
\end{defn}

\begin{figure}
    \centering
        \begin{subfigure}[b]{.49\textwidth}
       \includegraphics[scale=.45]{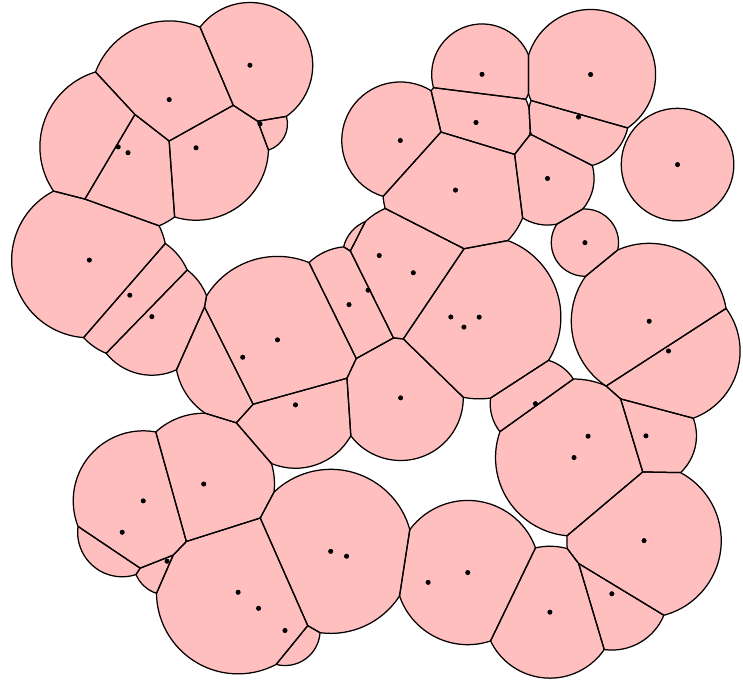}     
    \end{subfigure}
    \begin{subfigure}[b]{.49\textwidth}
       \includegraphics[scale=.45]{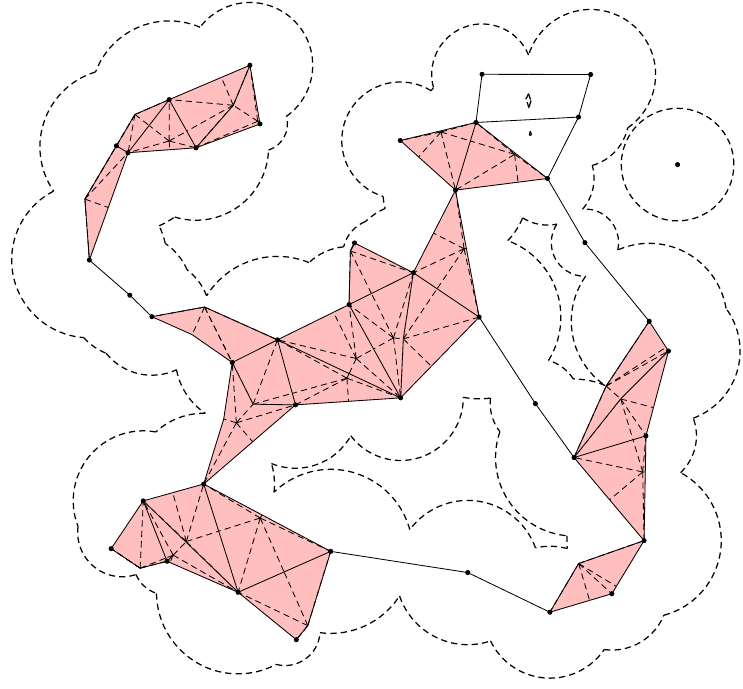}     
    \end{subfigure}
    \caption{On the left, a randomly generated power diagram in the plane, cut off at some weight $a$. On the right, the Barycentric subdivision of its associated alpha complex mapped into $\mathbb{R}^2$ using the witness map $\Phi$, which
    induces the nerve isomorphism.}
    \label{fig:nerve}
\end{figure}

In particular, by restricting by to $X(a)$, 
we obtain a linear map 
\[|\bary(X(a))|\rightarrow \bigcup \powdiag(S,p,a)\]
by the discussion in Section \ref{sec:nerve}.
An example of a power diagram, its alpha complex,
and the associated witness map is shown in Figure \ref{fig:nerve}.

The usual (unweighted) alpha complex
and Voronoi diagram from the introduction 
are given by
\[\alphaplex(S,r)=\alphaplex(S,p,r^2),\ 
\vor(S,r)=\powdiag(S,p,r^2)\]
in which $p(x)=0$ for all $x$. Notice that
the full vertex set for $X=\alphaplex(S,r)$ is given by
$X_0(a)=S$ for $a\geq 0$, and is empty for $a<0$, whereas the vertices appear at different times in the weighted case.
\section{Algorithm for computing the alpha complex}

We recall some facts about dual quadratic programming, and present our main algorithm.

\subsection{Dual programming}

\label{sec:dualprog}

Let $A$ be an $n\times m$ matrix, writing $A_i$ for the $i$th row. Let
$V\in \mathbb{R}^n$, $x_0\in \mathbb{R}^m$, and let $J\subset \{1,...,n\}$ be a subset.
Let $(y^*,c^*)=\primalqp(x,A,V,J)$
be the optimal solution and objective value to the quadratic program
\begin{mini}|l|
  {y\in \mathbb{R}^m}{\frac{1}{2}\lVert y-x\rVert^2}{}{}
  \addConstraint{A_{i} y}{=V_i,}{\quad i\in J}
  \addConstraint{A_{i} y}{\leq V_i,}{\quad i \notin J}
  \label{mini:primal}
 \end{mini}
To indicate that there is no feasible solution, the routine will return $c^*=\infty$, and an arbitrary value of $x^*$.

The dual quadratic program is
\begin{maxi}|l|
  {\lambda \in \mathbb{R}^n}{-\frac{1}{2}\lambda^t B\lambda+U^t \lambda}{}{}
  \addConstraint{\lambda_i{\geq 0,}{\quad i \notin J}}
  \label{maxi:dual}
 \end{maxi}
where
\begin{equation}
\label{eq:primtodual}
B=AA^{t},\quad U=A x-V.
\end{equation}
We will denote its solution by
$(\lambda^*,c^*)=\dualqp(B,U,J,c_1)$,
where $c_1$ is an upper bound on the allowable value of $c^*$. If we find that $c^*>c_1$, then the algorithm will terminate early and return $c^*=\infty$, together with an arbitrary value $\lambda^*$. If an optimum is obtained, the minimizing solution of the primal QP 
is determined by the KKT conditions, and is given by
\begin{equation}
    \label{eq:kkt}
y^*=x-A^t \lambda^*.    
\end{equation}

\begin{example}
\label{ex:lowerbound}
For every feasible point of \eqref{maxi:dual}, weak duality states that the corresponding value of the objective function is 
a lower bound on the solution in \eqref{mini:primal}.
Suppose for some $i$ we have that $A_{i}x< V_i$, meaning that $x$ is not on the feasible side in \eqref{mini:primal}. Then the minimizer of \eqref{maxi:dual} along the line $\lambda_i=t$ and all other $\lambda_j$ are zero occurs at
$t=U_i/B_{i,i}$. Substituting this into \eqref{maxi:dual} gives
\[-\frac{1}{2} B_{i,i}t^2+U_it=
\frac{U_i^2}{2B_{i,i}}=\frac{(A_{i}x-V)^2}{2\lVert A_{i}\rVert^2},\]
which is the lower bound corresponding to the point on the plane $A_{i}y=V$ which is as close as possible to $x$.
\end{example}

To solve \eqref{maxi:dual}, we have incorporated a compiled MAPLE implementation of a highly efficient and elegant recent active set method due to \cite{arnstrom2022dual}. 

\subsection{The alpha complex as a quadratic program}  

Consider the filtered alpha complex 
$X=\alphaplex(S,p,a_1)$ for $S\subset \mathbb{R}^m$, and let $w:X\rightarrow \mathbb{R}$ be the weight map.
Let $\sigma=[x_0,...,x_k] \subset S$ be a subset which may or may not define a simplex in $X$,
and and select any particular vertex, say $x=x_0$.
Then the problem of determining whether $\sigma$
determines a simplex in $X$ 
amounts to solving the following constrained quadratic optimization problem:
\begin{mini}|l|
  {y\in \mathbb{R}^m}{\pi_{x}(y)}{}{}
  \addConstraint{\pi_{x}(y)}{=\pi_{z}(y),}{\quad z\in \sigma-\{x\}}
  \addConstraint{\pi_{x}(y)}{\leq \pi_{z}(y),}{\quad z \in S-\sigma}
  \label{mini:pow}
 \end{mini}
Specifically, we have a simplex $\sigma \in X$ if and only if \eqref{mini:pow} is feasible and the optimal solution $a^*$ satisfies $a^*\leq a_1$. By convexity, if the problem is feasible, then there is a unique minimizer $y^*$ with corresponding objective value $a^*$, and we have $w(\sigma)=a^*$, and $\witmap(\sigma)=y^*$.

This can be formulated in terms of \eqref{mini:primal}. Let us write
$S-\{x\}=\{x_1,...,x_{n}\}$, and let $J=\{j_1,...,j_k\}$ be
those labels so that $\sigma-\{x\}=\{x_{j_1},...,x_{j_k}\}$. The contraints may be written as
\[
    A_{i}= (x_i-x)^t,\quad
    V_i= \frac{1}{2}\left(\lVert x_i\rVert^2-\lVert{x}\rVert^2-p(x_i)+p(x)\right)\]
for $i\in \{1,...,n\}$. Using \eqref{eq:primtodual}, the dual problem is determined by
\[B_{i,j}=(x_i-x)^t (x_j-x),\quad 
U_i=\frac{1}{2}\left(p(x_i)-p(x)-\lVert x_i-x\rVert^2\right).\]
Thus, if we set $(\lambda^*,c^*)=\dualqp(B,U,J,c_1)$ for $c_1=(a_1+p(x))/2$, then $\sigma$ determines a simplex if $c^*\leq c_1$, and its weight is given by
$w(\sigma)=2c^*-p(x)$.
The witness is given by $\Phi(\sigma)=y^*$, where by the KKT conditions \eqref{eq:kkt} we have
\[y^*=x-\sum_{i} \lambda^*_i(x_i-x).\]

\subsection{Description of the main algorithm}

We present our main algorithm which computes the 
weighted alpha complex using dual programming.
Specifically, the input consists of
 the data of a power diagram $(S,p,a_1)$ as in Section \ref{sec:powdiag}, together with a nonnegative integer $d\geq 0$.
 The output is the $d$-skeleton 
 $X=\skel_d(\alphaplex(S,p,a_1))$, which constains the simplices of the alpha complex up to dimension $d$, as well as 
 the associated witness map $\Phi:X\rightarrow \mathbb{R}^m$. 
 We now describe the procedure.

As a preprocessing step, we begin by computing the \v{C}ech 
graph $G$ of the weighted ball cover
$\powcov(S,p,a_1)=\{U_x:x\in S\}$, 
which is the one-skeleton of the \v{C}ech complex $\cechplex(S,p,a_1)$.
In other words, $G$ is the graph whose vertices are those elements $x\in S$ for which $U_x$ is nonempty, i.e. 
$-p(x)\leq a_1$, and which has an edge connecting $x$ to
$y$ if $U_x\cap U_y \neq \emptyset$. This is 
a important for efficiency for the following reason: 
let $x\in S$ and suppose $S'=\{x\} \cup N_G(x)\subset S$ 
contains $x$ and its neighbors in $G$. 
Then any face $V_\sigma$ of $V_x \in \powdiag(S,p,a_1)$ is nonempty if and only the corresponding face is nonempty in $\powdiag(S',p,a_1)$. This reduces the number of inequalities that need to be consider to only the necessary ones.

We then proceed to compute the alpha complex, beginning with the zero simplices, and work upwards until we reach dimension $d$. 
For each dimension $k\leq d$, we assume we have already computed the simplices of dimension up to $k-1$, described by a $(k-1)$-dimensional complex $X$, starting with the empty complex.
At step $k$, the algorithm may be described as follows.
\newpage
\begin{enumerate}
    \item Compute the set $\Sigma_k$ consisting of all potential simplices to be tested according to the following cases:
    \begin{enumerate}
        \item If $k=0$, then $\Sigma_k$ corresponds to
        the vertices of $G$. In other words, it is the set of simplices $[x]$ where $x\in S$ satisfies $-p(x)\leq a_1$.
\item If $k=1$, then $\Sigma_k$ is the set of edges of $G$.
\item If $k\geq 2$, then $\Sigma_k$ is the set of $k$-dimensional simplices in the lazy construction, $\Sigma_k=(\lazy_{k-1}(X))_k$. 
    \end{enumerate}
    \item For each $x\in S$, do the following:
    \begin{enumerate}
        \item Let $S'=\{x,x_1,...,x_n\}$ consist of $x$ together with its neighbors in $G$. Determine the coefficients $(B,U,c_1)$ 
        for the dual program that describes the cell $V_x \in \powdiag(S',p,a_1)$
        as in Section \ref{sec:dualprog}. The equations 
        defining every face of
        $V_x$ have the same coefficients, but 
        with different sets 
        $J$ that label the equality constraints.
        \item For every potential face $V_\sigma$ for $\sigma=[x,x_{j_1},...,x_{j_k}]\in \Sigma_k$, 
        do the following, noting that we are only considering
        simplices $x\leq x_{j_1}\leq \cdots \leq x_{j_k}$, as the other orders have already been encountered:
\begin{enumerate}
    \item \label{dual-active-set} Let $J=\{j_1,...,j_k\}$ be the indices of the
    equality constraints which determine the face $V_\sigma$.
    Solve the corresponding quadratic program using a dual active set method such as \cite{arnstrom2022dual}, terminating early if the upper bound $c_1$ is exceeded.
    \item If the quadratic program is feasible and the optimal value
    satisfies $c^*\leq c$,
    then add the simplex $\sigma$ to 
    $X_k$ with the desired weight $w(\sigma)$ 
    as determined by $c^*$. Compute the corresponding minimizer $x^*$ using the KKT equations, and update the witness map by setting $\Phi(\sigma)=x^*$.
\end{enumerate}
    \end{enumerate}
\end{enumerate}

\begin{figure}
\begin{centering}
\begin{subfigure}[b]{.24\textwidth}
\center
    \includegraphics[scale=.3]{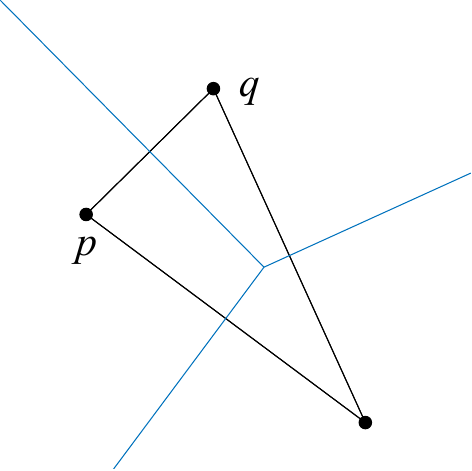}
        \caption{}
     \label{fig:3-points}
\end{subfigure}
\begin{subfigure}[b]{.24\textwidth}
\center
    \includegraphics[scale=.3]{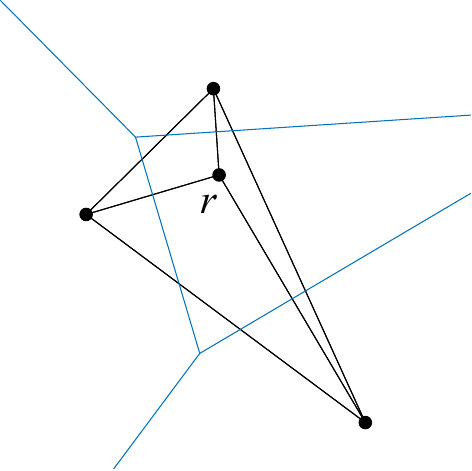}
        \caption{}
     \label{fig:4-points-1}
\end{subfigure}
\begin{subfigure}[b]{.24\textwidth}
\center
    \includegraphics[scale=.3]{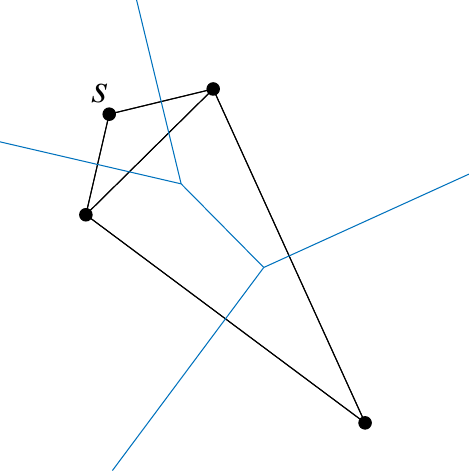}
    \caption{}
     \label{fig:4-points-2}
\end{subfigure}
\begin{subfigure}[b]{.24\textwidth}
\center
    \includegraphics[scale=.3]{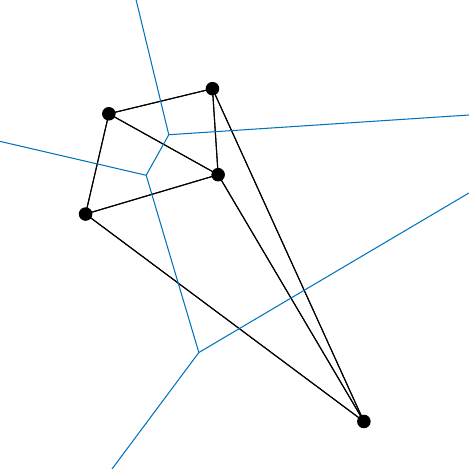}
        \caption{}
     \label{fig:5-points}
\end{subfigure}
\caption{The Delaunay triangulation (and Voronoi
diagram) for three points is shown in \ref{fig:3-points}. The segment
connecting points $p$ and $q$ continues to be present when points
$r$ and $s$ are inserted individually as in \ref{fig:4-points-1}
and \ref{fig:4-points-2}, but vanishes when both are inserted simultaneously
as in \ref{fig:5-points}. The key observation is that the simplex
$[p,q]$ can be ruled out due to the ``active set'' $\{r,s\}$,
irrespective of whatever other data points remain elsewhere in the point set.}
\label{fig:active-set-pic}
\end{centering}
\end{figure}

By using a dual active set method in step \ref{dual-active-set}, we can often rule out a potential simplex using a small subset of points, as shown in Figure \ref{fig:active-set-pic}, rather than all $N$ of them.  The process of solving problem (\ref{maxi:dual}) involves sequentially inserting and removing iterates $\lambda_i$, which are dual variables corresponding to the data points $x_i$, and this insertion and removal is equivalent to efficiently selecting a (typically small) subset of data points whose existence rules out a potential simplex. The psuedo-code for this algorithm is given in Algorithm \ref{alg:dualalpha}. 

\begin{algorithm}
    \caption{Compute the $d$-skeleton of $\alphaplex(S,p,a_1)$}
    \label{alg:dualalpha}    \begin{tabular}{lll} 
 \textbf{Input} &  $S\subset \mathbb{R}^m$ & Vertices of the power diagram\\ 
 &  $p:S\rightarrow \mathbb{R}$ & Power function\\ 
& $a_1\in\mathbb{R}$ & Maximum allowable power\\
& $d\geq 0$ & Dimension of output\\
\textbf{Output} & $X\subset \powset(S)$ & $d$-skeleton of $\alphaplex(S,p,a_1)$\\
& $w:X\rightarrow (-\infty,a_1]$ & Weight function \\
& $\witmap : X \rightarrow \mathbb{R}^m$ & Table of representatives 
\end{tabular}
    \begin{algorithmic}[1]
\State $Y\gets \skel_1(\cechplex(S,p,a_1))$
\State $G \gets \graph(Y)$
\algorithmiccomment{Underlying graph}
\label{line:graph}
\State $X \gets \emptyset$
\For{$k$ from 0 to $d$}
    \If{$k\leq 1$} 
    \State $\Sigma_k \gets Y_k$        
    \Else
           \State$\Sigma_k \gets (\lazy_{k-1}(X))_k$ \algorithmiccomment{All potential $k$-simplices}
    \EndIf
    \For{$x \in S$}
         \State $\{x_1,...,x_n\} \gets \neighbors_G(x)$ \algorithmiccomment{Neighbors to $x$}   
        \State \label{line:setb} $B \gets ((x_{i}-x)^t(x_j-x))_{i,j=1}^{n}$  \label{line:coeffs}
        \State \label{line:setu} $U\gets \frac{1}{2}(p(x_i)-p(x)-\lVert x_{i}-x \rVert^2)_{i=1}^{n}$    
        \State $c_1\gets (a_1+p(x))/2$
        \For{$\sigma=[x,x_{j_1}...,x_{j_k}] \in \Sigma_k$}
         \State $J\gets \{j_{1},...,j_{k}\}$ 
         \algorithmiccomment{The equality constraints}
    \State $(\lambda^*,c^*) \gets \dualqp(B,U,J,c_1)$ 
    \If{$c^*\leq c_1$}
        \State $X \gets X \cup \{\sigma\}$
        \State $w(\sigma)\gets 2c^*-p(x)$
        \State $\witmap(\sigma)\gets x-\sum_{j=1}^{n} \lambda^*_j (x_j-x)$
        \algorithmiccomment{KKT conditions}
    \EndIf
            \EndFor 
    \EndFor
\EndFor
\end{algorithmic}
\end{algorithm}

We summarize the above discussion in a proposition, 
which is evident:
\begin{prop}
\label{prop:alg}    
Algorithm \ref{alg:dualalpha} computes the alpha complex.
\end{prop}

\section{Examples and Applications}

\label{sec:examples}

We illustrate Algoritithm \ref{alg:dualalpha} in several examples. In the first, we apply the complex to a standard three-dimensional mesh generation example. We find that Algorithm \ref{alg:dualalpha} is not as fast as 
existing methods that are specialized to three dimensions. 
In the second, we generate 1000 random points in $\mathbb{R}^{10}$ and compute the alpha complex up to the three-dimensional simplices with a relatively large radius. In contrast with the three-dimensional example, this could not be done by computing the full Delaunay triangulation.

In the last two examples, we use the alpha complex to compute the 
homology of some interesting topological spaces 
from a sampling of landmark points. 
While the main loop of Algorithm \ref{alg:dualalpha} can be done in full parallel, we have not used any parallelism in our computations.
The homology groups calculated below were done using a MAPLE implementation of an algorithm of Dumas and Villard \cite{dumas2002sparse} for computing the ranks of sparse matrices mod $p$. For comparison, we also compute the persistence homology groups using Ripser.

\subsection{Three-dimensional mesh generation}

\begin{figure}
    \centering
    \includegraphics[scale=.3]{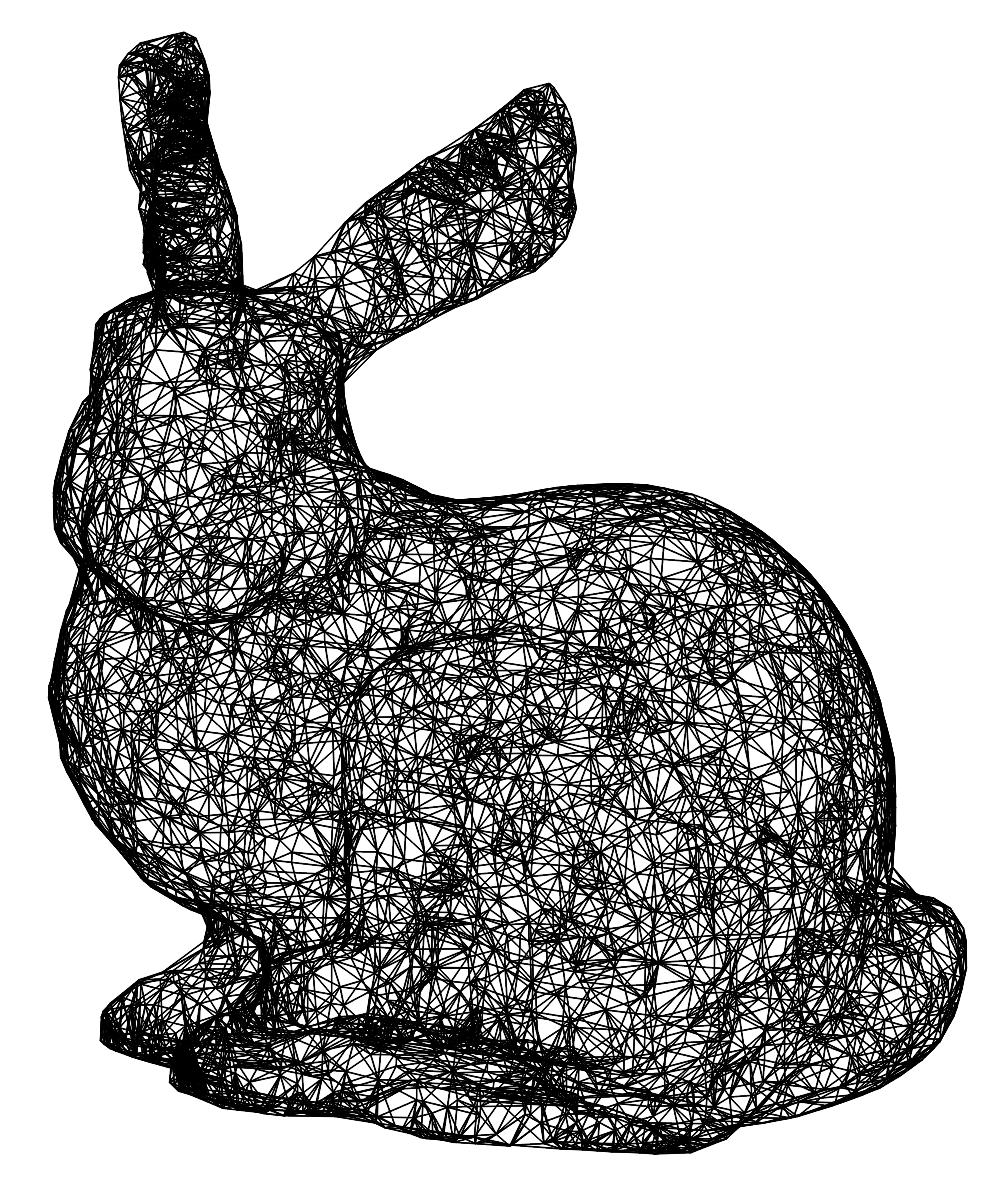}
    \caption{Alpha complex of a Stanford bunny data set with 5000 sites.}
    \label{fig:bunny}    
\end{figure}

We first apply Algorithm \ref{alg:dualalpha} to triangulating 
a three dimensional data set consisting of
 5000 points sampled from surface of the Stanford bunny \cite{turk1994zippered}, downloaded from the CGAL website \cite{cgal:eb-23a}. We chose a radius size of approximately 1/15 of the diameter, which led to a \v{C}ech graph with a maximum vertex degree of 200 in the graph of line \ref{line:graph}. The resulting alpha
complex, whose one-skeleton is shown in Figure \ref{fig:bunny}, had sizes of $(|X_k|)_{k=0}^3=(5000,23830,29995,11163)$ simplices in each dimension. The full computation took approximately 12 seconds. This would be faster to compute using either specialized methods for three-dimensions, or using Delaunay software such as qhull \cite{qhull}. Notice that the Euler characteristic gives the value of 2, indicating that the corresponding covering is homotopy equivalent to the sphere.

\subsection{Random points in higher dimension}

The next example could not be done by computing the full Delaunay triangulation. We chose 1000 points in $\mathbb{R}^{10}$
by selecting each coordinate uniformly at random from the interval $[-1,1]$, producing a point cloud $S\subset \mathbb{R}^{10}$. We then computed the alpha complex $X=\alphaplex(S,1.0)$ up to the 3-simplices.
The radius was such that the maximum degree in the \v{C}ech graph was about half the size of the data set, in our example 451. The full calculation took approximately 5 minutes, resulting
in complex of sizes $(|X_k|)_{k=0}^3=(1000,64785,560459,1783194)$.

\subsection{Two persistence examples}

We next tested our algorithm on two well-studied data sets from persistent homology. Both are freely available online, and are explained in Henry Adams' tutorial on topological data analysis and
Ripser \cite{adams2000wiki}. We find that the first one, which is a 24-dimensional data set consisting of conformations of the cyclooctane molecule, is well-suited for the alpha complex because it tends to lie near the surface of a lower-dimensional space. The second one, which is a $9$-dimensional database of optical image patches, leads to a complex with more simplices, because it has thickness in more dimensions, despite being embedded in lower dimensions.

\begin{figure}
    \centering
    \includegraphics[scale=.5]{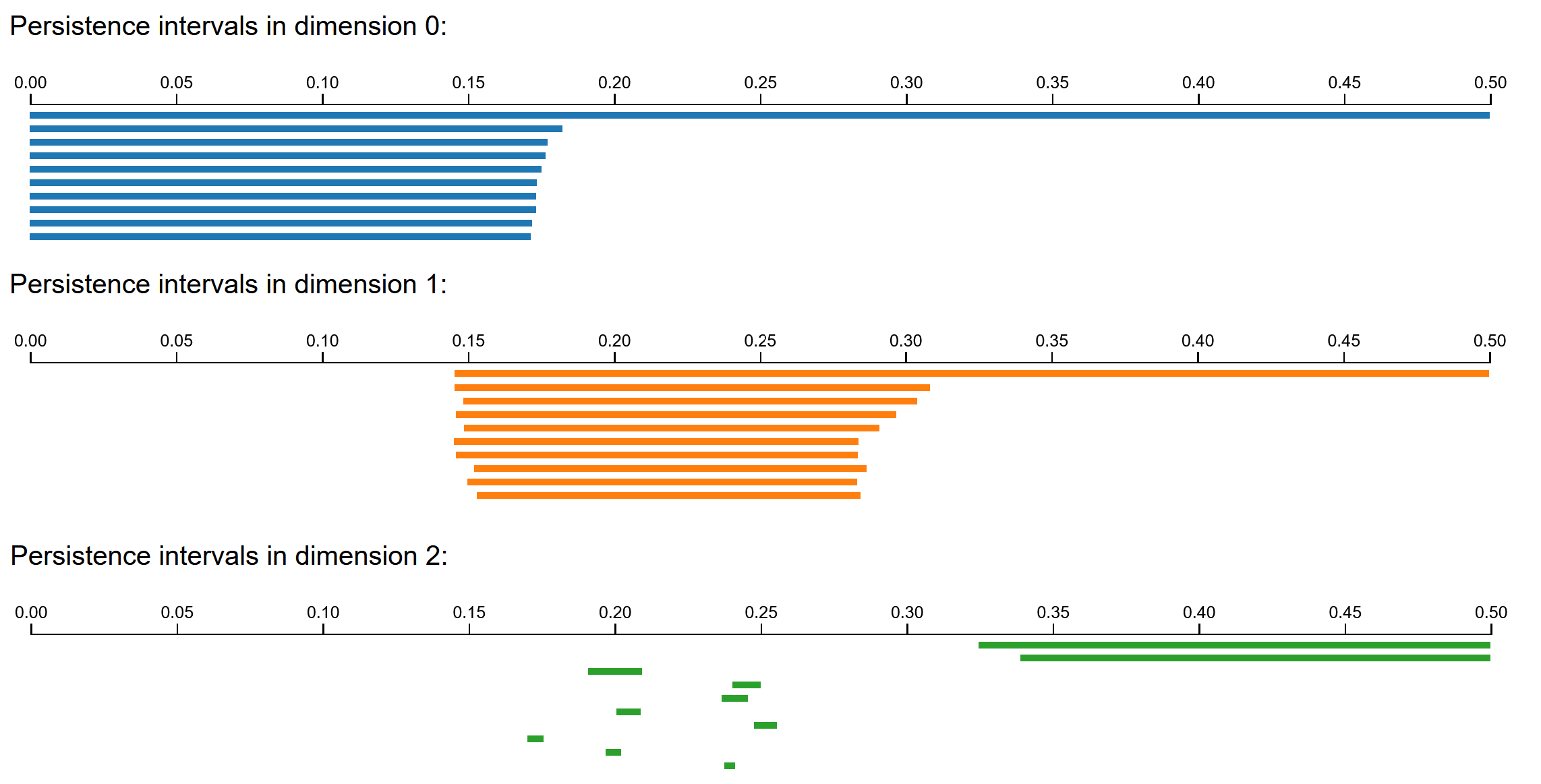}
    \caption{Persistence intervals of the cyclooctane data set as computed by Ripser. Many intervals not shown for viewability.}
    \label{fig:cyclorips}
\end{figure}

Our first example is a data set consisting of 6040 points $S\subset \mathbb{R}^{24}$ in 24 dimensions, which correspond to conformations of the cyclooctane molecule, introduced in \cite{martin2010topology}. The authors found that the set of conformations 
tend to lie on a 2-dimensional topological subspace $X$ which is an interesting union of a Klein bottle and a sphere. Because of its interesting topological type, it is a well-suited use case of persistent homology, in particular Ripser. By running Ripser on the full data set up to a cutoff distance of $.5$, we obtain a diagram of persistence intervals as shown in Figure \ref{fig:cyclorips}, which agrees with the desired Betti numbers of $(\betti_k(X))_{k=0}^2=(1,1,2)$. Ripser took approximately 48 seconds to complete this calculation.

We then computed the alpha complex
$\alphaplex(S,.25)$ up the $3$-simplices, as required to compute up to the second Betti number.
Notice that our cutoff of $a_1=.25$ is half the cutoff used for Ripser, because the minimum radius at which two balls intersect is half the distance between the centers. This computation took approximately 7 seconds, producing a complex of sizes
$(|X_k|)_{k=0}^3=(6040,24646,28352,11858)$.
We checked that it produced the desired Betti numbers exactly, without any persistence.
Despite the much smaller size of the alpha complex as compared with the Vietoris-Rips construction, our homology calculation took longer, as we made no effort to use specialized algorithms.

We then applied the same procedure to a nine-dimensional data set studied in \cite{lee2003nonlinear}, consisting of normalized $3\times 3$ patches taken from the van Hateren and van er Shaaf image database \cite{vanhateren1998independent}. Certain high-density subsets were studied using persistent homology in \cite{carlsson2008klein}, revealing the topological type of certain subspaces of a parametrized Klein bottle. We apply our algorithm to a subset $S\subset \mathbb{R}^9$ of size 1000 which is known to have the topological type of a circle, and which is
denoted $X(300,30)$ in \cite{carlsson2008klein}.
We computed $\alphaplex(S,.5)$ up to the 2-dimensional simplices, obtaining a complex of sizes
$(|X_k|)_{k=0}^2=(1000,44080,454843)$ in about 38 seconds, and computed Betti numbers of $(1,1)$, which agree with those of the circle. In this example, Ripser took only 1.5 seconds to obtain persistence intervals indicating these numbers. 

\subsection{Spherical images from different angles}

We next consider the alpha complex of a data set consisting of $28\times 28$ color images of a coloring of the surface of the sphere from different angles, viewed as vectors in dimension $28\times28\times 3$. 
This example could not be accomplished 
using an algorithm that 
begins by computing the full Delaunay triangulation
due to the prohibitively high dimension.
This illustrates a point that in dual programming the high dimensionality of the ambient space is not a direct factor as the input is a function only of respective dot products. Indeed, the only part of Algorithm \ref{alg:dualalpha} that depends explicitly on the dimension is line \ref{line:coeffs}, in which the coefficients of the quadratic program are computed, 
which only happens once per vertex in each dimension.

\begin{figure}
    \centering
    \begin{subfigure}[b]{.24\textwidth}
         \includegraphics[scale=1.0]{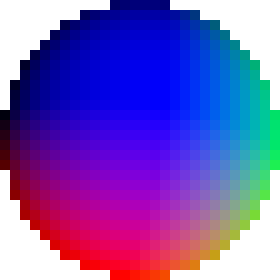}   
    \end{subfigure}
       \begin{subfigure}[b]{.24\textwidth}
         \includegraphics[scale=1.0]{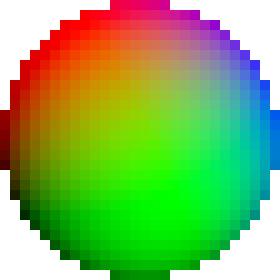}   
    \end{subfigure}
       \begin{subfigure}[b]{.24\textwidth}
         \includegraphics[scale=1.0]{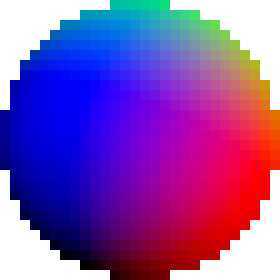}   
    \end{subfigure}
       \begin{subfigure}[b]{.24\textwidth}
         \includegraphics[scale=1.0]{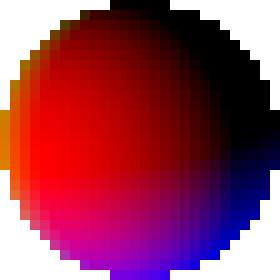}   
    \end{subfigure}

\vspace{.1in}
    
        \begin{subfigure}[b]{.24\textwidth}
         \includegraphics[scale=1.0]{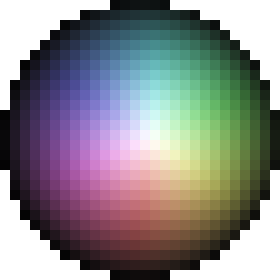}   
    \end{subfigure}
       \begin{subfigure}[b]{.24\textwidth}
         \includegraphics[scale=1.0]{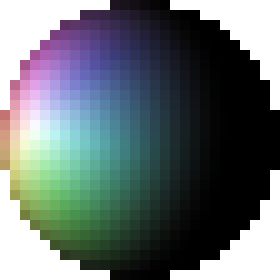}   
    \end{subfigure}
       \begin{subfigure}[b]{.24\textwidth}
         \includegraphics[scale=1.0]{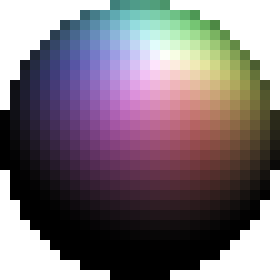}   
    \end{subfigure}
       \begin{subfigure}[b]{.24\textwidth}
         \includegraphics[scale=1.0]{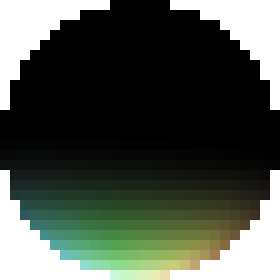}   
    \end{subfigure}

    \caption{Top row: Pictures of a sphere from colored by $\varphi_1$ from different angles, corresponding the RGB values with the three axes. Bottom row: The same for $\varphi_2$, displayed using a cylindrical HSV coloring scheme.}
    \label{fig:so3}
\end{figure}

Fix a function $\varphi : S^2\rightarrow \mathbb{R}^3$ thought of as a coloring of the surface of a sphere. We will be interested in the following two choices:
\[ \varphi_1(x,y,z)=(x,y,z),\quad \varphi_2(x,y,z)=(tx,ty,tz)\]
where $t=\max(z,0)$, so that one hemisphere is sent to the origin.
We then have a map $F_{\varphi}:SO(3)\rightarrow \mathbb{R}^m$
for $m=28\cdot 28\cdot 3$
defined by projecting $\varphi\circ R^{-1}$ onto the $xy$-plane
and discreting the result into a $28\times 28$ image. In other words, we take an image using a camera with position defined by $R$ with no perspective warping.
A collection of these images are shown in Figure \ref{fig:so3}.

For $\varphi=\varphi_1$, we randomly generated $N=1000000$ images, producing a point cloud $D\subset \mathbb{R}^m$. We then selected landmarks 
$S\subset X$, until every element of $D$ was within a distance of 8.0 of some element of $S$, giving 
a size of $|S|=809$.
Assuming that a covering by balls of radius $10.0$ with centers at the points of $S$ would cover the entire image of $F_{\varphi_1}$, we proceeded to compute $X=\alphaplex(S,10.0)$ up to the $4$-simplices. The resulting complex took approximately 20 seconds in each dimension, and had sizes of
\[(|X_k|)_{k=0}^4=(809,7717,17694,15490,5599),\]
We then computed the Betti numbers up to dimension 3 over the fields $\mathbb{F}_2,\mathbb{F}_3$, giving
\[(\betti_k(X,\mathbb{F}_2))=(1,1,1,1),\quad (\betti_k(X,\mathbb{F}_3))=(1,0,0,1).\]
These are the desired Betti numbers of $SO(3)\cong \mathbb{RP}^3$, which is to be expected if $F_{\varphi_1}$ is a reasonable embedding of 
$SO(3)$ in $\mathbb{R}^m$.

We then performed a similar calculation for $\varphi=\varphi_2$.
Since distances are smaller for this embedding, we used a minimum distance of $3.0$ in selecting the landmark points, yielding a set $S$ of size
$2361$. We then computed $\alphaplex(S,3.5)$ up to the 4-simplices,
which took approximately fifteen minutes, 
yielding a complex of $X$ with sizes
\[(|X_k|)_{k=0}^4=(2340,24463,68150,102772,128302).\]
We first notice that the number of simplices increases more rapidly with degree, whereas the one from the previous paragraph began to decrease above the $2$-simplices. This is because the image $F_{\varphi_2}(SO(3))$ is pinched at the angle which is entirely on the dark side, resulting in a higher-dimensional tangent space. This makes the alpha complex effectively four-dimensional nearby the singularity. For instance, a single vertex corresponding to a nearly entirely black image had $(1,105,1756,11776,39376)$ simplices in each degree containing it as a face. 

The Betti numbers were the calculated as above, yielding
\[(\betti_k(X,\mathbb{F}_p))=(1,0,1,1)\]
for all primes $p$.
In fact, these are the desired Betti numbers resulting from collapsing 
the circle $S^1\subset SO(3)$ to a point, 
where the circle is the stabilizer of 
the image that is completely on the dark side of the coloring.
This agrees with what may be computed using the long exact sequence for relative homology groups $H_k(SO(3),S^1)$.
In the case of $p=2$, we must use the fact that the map 
$H_1(S^1,\mathbb{F}_2)\rightarrow H_1(SO(3),\mathbb{F}_2)$ is nonzero, 
so that the connecting homomorphism is trivial.
This results from the fact that $S^1$ corresponds to
a generator of the fundamental group $\pi_1(SO(3))=\mathbb{Z}_2$.

\begin{figure}
\includegraphics[scale=.5]{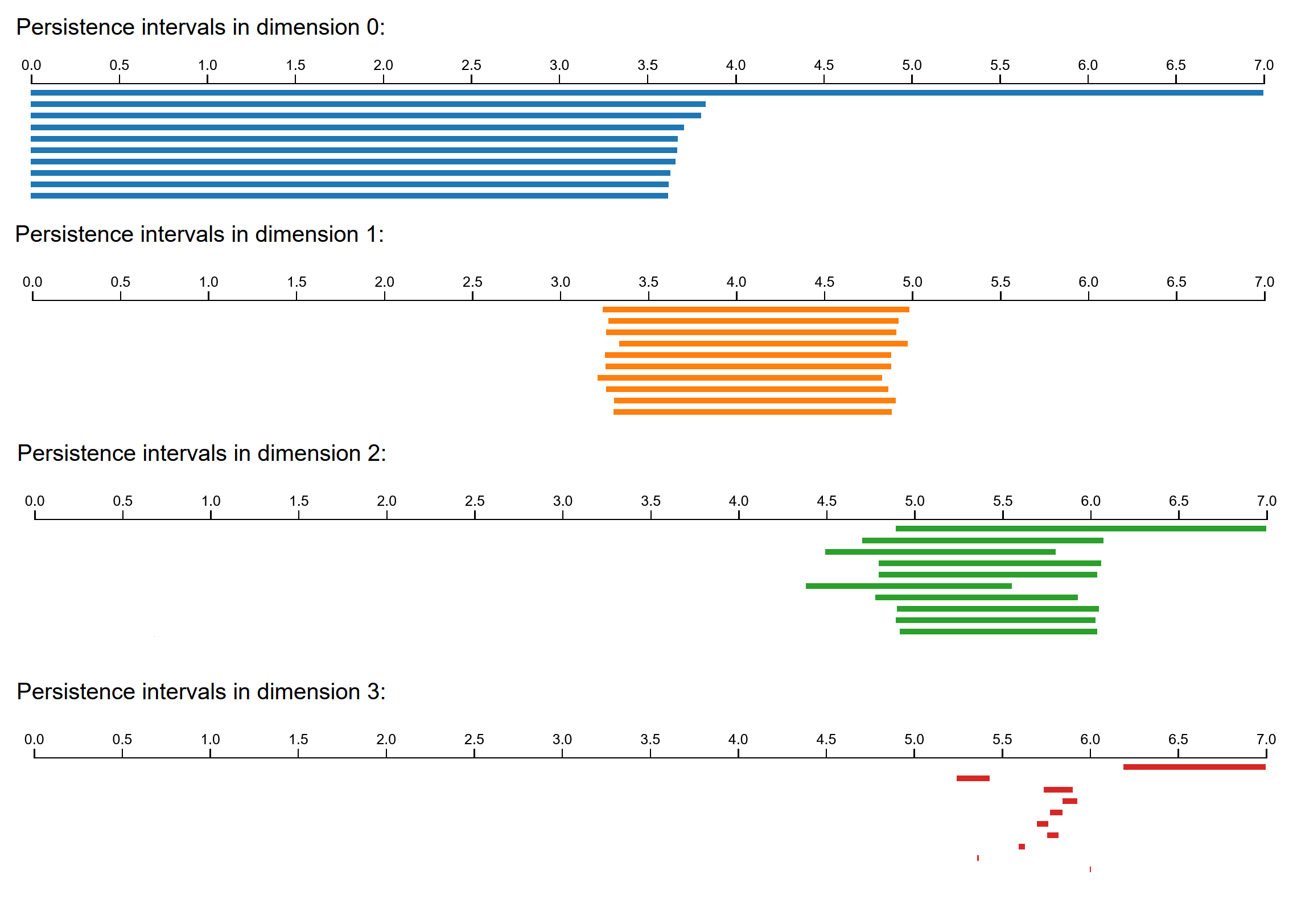}
\caption{Results of Applying Ripser to samples from the data set of spherical images with one side colored entirely black. The intervals which extend to the rightmost endpoint indicate Betti numbers of $(1,0,1,1)$.}
\label{fig:so3rips}
\end{figure}

We then uploaded the second example into Ripser up to a distance cutoff of 7.0. The results, which took approximately 30 seconds to compute, are shown in Figure \ref{fig:so3rips}.

\subsection{Homology groups of configuration spaces}

\label{sec:config}

In the next example, we use the alpha complex to carry out
a complex homology calculation, 
which is the homology groups of a compact form of the ordered
configuration space $\conf_n(\mathbb{R}^2)$ of $n$-points in the plane for $n=3,4$. Additionally, $\conf_n(\mathbb{R}^2)$ is 
acted on freely by the symmetric group $S_n$ by relabeling, 
which induces an action on homology.
By selecting landmark points in symmetric way, we generalize the Betti number computation and produce the decomposition into irreducible characters 
of the corresponding group representation, which is an interesting object in geometric representation theory. 
This makes use of the theoretically sound nature of the alpha complex, 
which computes homology exactly.

Let $\conf_n(M)$ denote the configuration space of $n$ points in $M$
\[\conf_n(M)=\left\{ (x_1,...,x_n) \in M^n: 
\mbox{$x_i\neq x_j$ for all $i\neq j$}\right\}.\]
Its homology and cohomology have been well-studied, see \cite{cohen1995configuration}.
Of particular interest is in the case of $\conf_n(\mathbb{R}^2)$, in which 
case the cohomology groups are the same as the
cohomology of the pure Artin Braid group
\cite{vainshtein1978cohomology}.
The Betti numbers are given 
in characteristic zero in this case by the formula
\begin{equation}
\label{eq:sterling} 
\sum_{k\geq 0} (-t)^k \betti_k(\conf_n(\mathbb{R}^2),\mathbb{Q})=
(1-t)(1-2t)\cdots (1-(n-1)t).
\end{equation}
For instance, the first three Betti numbers of $\conf_3(\mathbb{R}^2)$ would be $(1,3,2)$ with all higher Betti numbers being zero, and would be $(1,6,11,6)$ in the case of $n=4$. These numbers are known as Stirling numbers of the first kind. 

We also the action of the symmetric group by reordering the labels of the points,
\[\sigma\cdot (x_1,...,x_n)=(x_{\sigma_1},...,x_{\sigma_n}).\]
The induced action makes both homology 
and cohomology into representations of the 
symmetric group. They turn out to be graded versions of the regular representation
after twisting by the sign representation in odd degree, noticing that the above total dimensions from the previous paragraph sum to $n!$. A formula for the character of this action over complex coefficients is a special case of a results of Lehrer and Solomon \cite{lehrer1986action}. 
In this case, it says that the trace of the action of a permutation 
$\sigma$ on cohomology is given by
\begin{equation}
    \label{eq:lehrersolomon}
\sum_{i} (-t)^i \Tr(\sigma,H^i(\conf_n(\mathbb{R}^2),\mathbb{C}))=
t^n \prod_{i=1}^{\infty} \prod_{j=1}^{m_i(\lambda)} (\alpha_i(t^{-1})-(j-1)i)
\end{equation}
where $\lambda=(\lambda_1,...,\lambda_l)$ is a Young diagram 
describing the cycle type of $\sigma$, $m_i(\lambda)$ is the number of times that each $i$ appears in the elements of $\lambda$, and 
\[\alpha_j(t)=\sum_{d|j} t^d\mu(j/d)\]
where $\mu$ is the usual M\"{o}bius function. In particular, we recover formula \eqref{eq:sterling} by taking $\sigma$ to be the identity permutation.

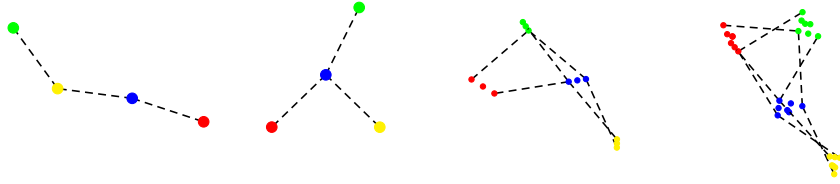
\begin{figure}
    \centering
    \begin{subfigure}[b]{.24\textwidth}
    \centering
    \begin{tikzpicture}[scale=1.0]
    \draw[color=white] (-1.3,-1.3)--(1.3,1.3);
\draw[black,semithick,densely dashed](0.405534,-0.188474)--(1.355466,-0.500930);
\draw[black,semithick,densely dashed](0.405534,-0.188474)--(-0.586107,-0.059442);
\draw[black,semithick,densely dashed](-1.174893,0.748847)--(-0.586107,-0.059442);
\filldraw[blue] (0.405534,-0.188474) circle (2.000000pt);
\filldraw[red] (1.355466,-0.500930) circle (2.000000pt);
\filldraw[green] (-1.174893,0.748847) circle (2.000000pt);
\filldraw[yellow] (-0.586107,-0.059442) circle (2.000000pt);
\end{tikzpicture}        
    \end{subfigure}
    \begin{subfigure}[b]{.24\textwidth}
    \centering
\begin{tikzpicture}[scale=1.0]
    \draw[color=white] (-1.3,-1.3)--(1.3,1.3);
\draw[black,semithick,densely dashed](-0.110715,0.123948)--(-0.829265,-0.571527);
\draw[black,semithick,densely dashed](-0.110715,0.123948)--(0.333084,1.020075);
\draw[black,semithick,densely dashed](-0.110715,0.123948)--(0.606897,-0.572496);
\filldraw[blue] (-0.110715,0.123948) circle (2.000000pt);
\filldraw[red] (-0.829265,-0.571527) circle (2.000000pt);
\filldraw[green] (0.333084,1.020075) circle (2.000000pt);
\filldraw[yellow] (0.606897,-0.572496) circle (2.000000pt);
\end{tikzpicture}        
    \end{subfigure}
    \begin{subfigure}
        [b]{.24\textwidth}
        \begin{tikzpicture}
           \draw[color=white] (-1.3,-1.3)--(1.3,1.3);
\draw[black,semithick,densely dashed](0.468596,0.067800)--(-0.294652,0.713906);
\draw[black,semithick,densely dashed](0.468596,0.067800)--(0.879252,-0.843990);
\draw[black,semithick,densely dashed](-1.053196,0.062284)--(-0.294652,0.713906);
\filldraw[blue] (0.468596,0.067800) circle (1.000000pt);
\filldraw[red] (-1.053196,0.062284) circle (1.000000pt);
\filldraw[green] (-0.294652,0.713906) circle (1.000000pt);
\filldraw[yellow] (0.879252,-0.843990) circle (1.000000pt);
\draw[black,semithick,densely dashed](0.238843,0.032699)--(-0.748806,-0.123983);
\draw[black,semithick,densely dashed](0.238843,0.032699)--(-0.370534,0.825579);
\draw[black,semithick,densely dashed](0.238843,0.032699)--(0.880497,-0.734295);
\filldraw[blue] (0.238843,0.032699) circle (1.000000pt);
\filldraw[red] (-0.748806,-0.123983) circle (1.000000pt);
\filldraw[green] (-0.370534,0.825579) circle (1.000000pt);
\filldraw[yellow] (0.880497,-0.734295) circle (1.000000pt);
\filldraw[blue] (0.353719,0.050249) circle (1.000000pt);
\filldraw[red] (-0.901001,-0.030850) circle (1.000000pt);
\filldraw[green] (-0.332593,0.769743) circle (1.000000pt);
\filldraw[yellow] (0.879875,-0.789143) circle (1.000000pt);
\end{tikzpicture}
    \end{subfigure}
\begin{subfigure}[b]{.24\textwidth}
\centering
\begin{tikzpicture}
   \draw[color=white] (-1.3,-1.3)--(1.3,1.3);
\draw[black,semithick,densely dashed](0.169481,-0.292156)--(0.117152,0.706474);
\draw[black,semithick,densely dashed](0.169481,-0.292156)--(0.593235,-1.197933);
\draw[black,semithick,densely dashed](-0.879868,0.783615)--(0.117152,0.706474);
\filldraw[blue] (0.169481,-0.292156) circle (1.000000pt);
\filldraw[red] (-0.879868,0.783615) circle (1.000000pt);
\filldraw[green] (0.117152,0.706474) circle (1.000000pt);
\filldraw[yellow] (0.593235,-1.197933) circle (1.000000pt);
\draw[black,semithick,densely dashed](-0.158088,-0.415312)--(-0.681587,0.436714);
\draw[black,semithick,densely dashed](-0.158088,-0.415312)--(0.667759,-0.979207);
\draw[black,semithick,densely dashed](-0.681587,0.436714)--(0.171915,0.957804);
\filldraw[blue] (-0.158088,-0.415312) circle (1.000000pt);
\filldraw[red] (-0.681587,0.436714) circle (1.000000pt);
\filldraw[green] (0.171915,0.957804) circle (1.000000pt);
\filldraw[yellow] (0.667759,-0.979207) circle (1.000000pt);
\filldraw[blue] (-0.006742,-0.373150) circle (1.000000pt);
\filldraw[red] (-0.758263,0.635928) circle (1.000000pt);
\filldraw[green] (0.159275,0.845762) circle (1.000000pt);
\filldraw[yellow] (0.605730,-1.108540) circle (1.000000pt);
\draw[black,semithick,densely dashed](-0.135274,-0.220265)--(-0.778283,0.545593);
\draw[black,semithick,densely dashed](-0.135274,-0.220265)--(0.378939,0.637397);
\draw[black,semithick,densely dashed](-0.135274,-0.220265)--(0.534618,-0.962724);
\filldraw[blue] (-0.135274,-0.220265) circle (1.000000pt);
\filldraw[red] (-0.778283,0.545593) circle (1.000000pt);
\filldraw[green] (0.378939,0.637397) circle (1.000000pt);
\filldraw[yellow] (0.534618,-0.962724) circle (1.000000pt);
\filldraw[blue] (0.017104,-0.256211) circle (1.000000pt);
\filldraw[red] (-0.829076,0.664604) circle (1.000000pt);
\filldraw[green] (0.248046,0.671936) circle (1.000000pt);
\filldraw[yellow] (0.563927,-1.080329) circle (1.000000pt);
\filldraw[blue] (-0.146681,-0.317789) circle (1.000000pt);
\filldraw[red] (-0.729935,0.491153) circle (1.000000pt);
\filldraw[green] (0.275427,0.797601) circle (1.000000pt);
\filldraw[yellow] (0.601189,-0.970965) circle (1.000000pt);
\filldraw[blue] (-0.030633,-0.346087) circle (1.000000pt);
\filldraw[red] (-0.758852,0.632291) circle (1.000000pt);
\filldraw[green] (0.207422,0.803238) circle (1.000000pt);
\filldraw[yellow] (0.582063,-1.089442) circle (1.000000pt);
\end{tikzpicture}
\end{subfigure}    
    \caption{On the left, two typical points in $C_4$, one with each of the two types of associated tree structure describing which points are separated by distance exactly 1. On the right, the endpoints of the Barycentric subdivision of simplices of the resulting alpha complex $X^{(4)}$ of dimensions 1 and 2, embedded using the witness map $\Phi$. Notice that simplices may connect vertices with different associated graphs $G$.}
    \label{fig:config}
\end{figure}

A priori, $\conf_n(\mathbb{R}^2)$ is not well-suited to being triangulated because the constraint that $x_i\neq x_j$ is of measure zero, and so would not be detected by distance measurements.
We will instead replace the full configuration space with the compact subspace $C_n\subset \conf_n(\mathbb{R}^2)$, which consists of all points $(x_1,...,x_n)\in \mathbb{R}^2$ satisfying:
\begin{itemize}
    \item The points are mean centered, i.e. $x_1+\cdots +x_n=(0,0)$.
    \item For every $i\neq j$, we have that $\lVert x_i-x_j\rVert\geq 1$.
\item The graph $G$ on $n$ vertices which contains an edge connecting 
$i$ and $j$ whenever we have equality $\lVert x_i-x_j\rVert=1$ is connected.
\end{itemize}
We illustrate some typical points in Figure \ref{fig:config}, showing the graph $G$ in dashed lines, which is generically a tree.
While $\conf_n(\mathbb{R}^2)$ is smooth but not compact of dimension $2n$, it is not hard to see that $C_n$ is singular but compact of dimension $n-1$.
We also observe that $C_n$ is preserved by the $S_n$ action on $\conf_n(M)$.

In a similar way as previous example, we selected landmark points $S^{(n)}$ from $C_n$ for the values of $n=3,4$, using a distance cutoff of $.3$ in both cases. 
In order to maintain an $S_n$-action on the complex itself, we sampled in such a way that whenever a single point is added to, we also added its $S_n$-orbit. The resulting sizes were $|S^{(3)}|=328,|S^{(4)}|=20232$, noticing 
that the sizes are multiples of 6 and 24 respectively.

\begin{figure}
\includegraphics[scale=.5]{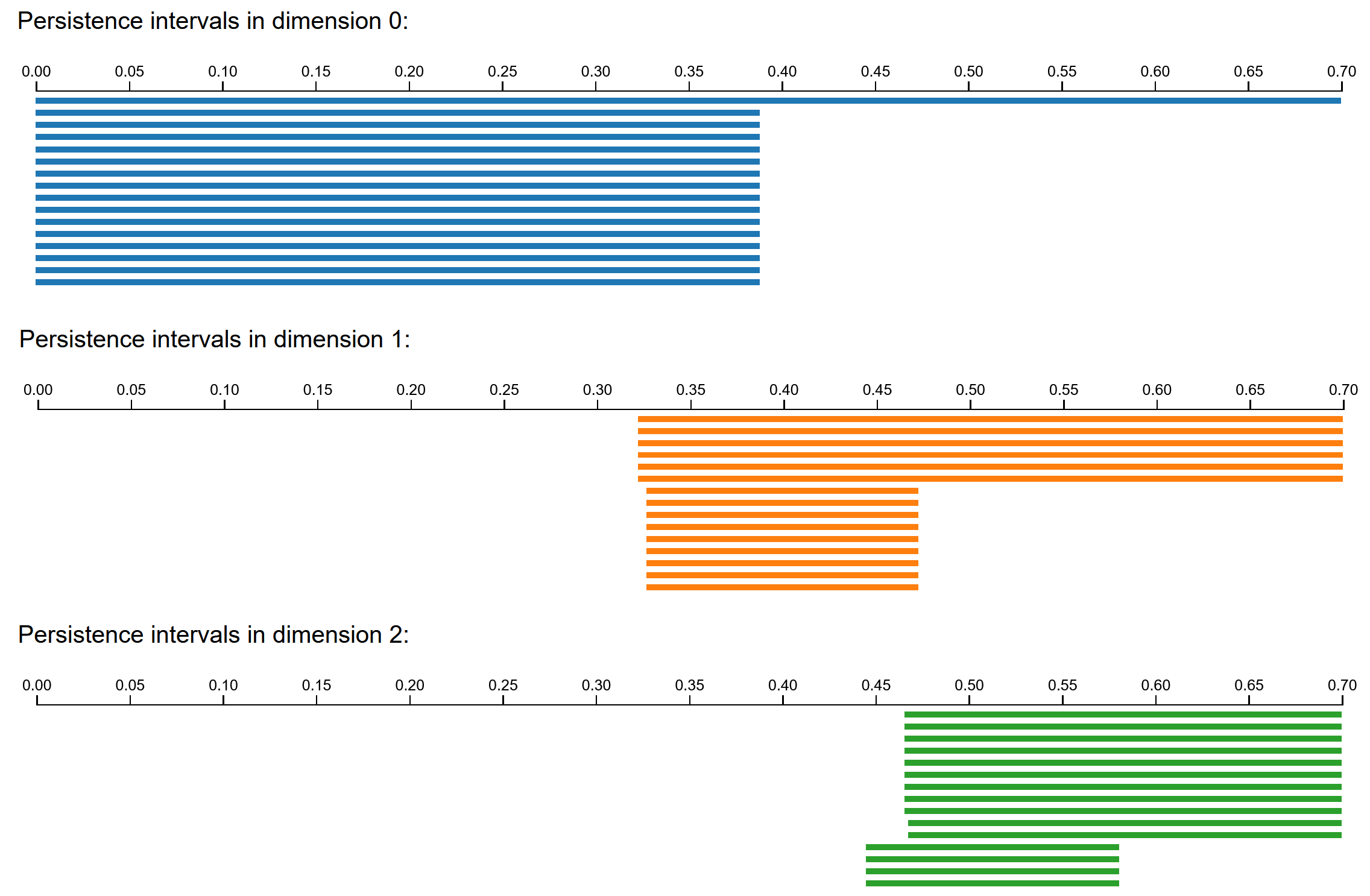}
\caption{Results of Applying Ripser to samples from the configuration space $C_4$, with a very large number of shorter length persistence intervals excluded for viewability. We see the first three Betti numbers of $(1,6,11)$.}
\label{fig:configrips}
\end{figure}

We first attempted to compute the homology groups using Ripser, by uploading $S^{(n)}$ to the Ripser live, with a cutoff distance of $.7$, which more than the lower bound of .6 for the distance between any two points. Not surprisingly, the case of $n=3$ was trivial. Ripser was able to compute the homology groups up to the second Betti number for $S^{(4)}$ in 30-40 minutes, though did not manage to compute the third Betti number. The results are shown in Figure \ref{fig:configrips}.

We then computed the full alpha complexes to top dimension, with a distance cutoff of .35, yielding filtered complexes $X^{(n)}=\alphaplex(S^{(n)},.35,4)$. The first one $X^{(3)}$ took under one second to compute, whereas computing $X^{(4)}$ took approximately 4 minutes. 
 The resulting sizes were
 \begin{equation} \label{eq1}
\begin{split}
(|X^{(3)}|)_{k=0}^4 & = (328,1245,1201,312,28),\\
 (|X^{(4)}|)_{k=0}^4 & = (20232,228072,608928,665208,345888,92040,10272),
\end{split}
\end{equation}
all others sizes being empty because the spaces are embedded into dimension 
$2n-2$ by the mean-centering condition.
Notice that the Euler characteristics are zero, and that each number is a multiple of $n!$.
Some typical two and three-simplices are shown on the right side of Figure \ref{fig:config}.
By the way we chose the landmark points, we have an action of the symmetric group on each $X^{(n)}_k$, and that the boundary maps commute with this action.

The two spaces had Betti numbers that agreed with \eqref{eq:sterling}.
To compute the character of the $S_n$-representation, 
we applied Mashke's theorem to convert the boundary operator $\partial_k : C_k(X^{(n)},\mathbb{F}_5)\rightarrow C_{k-1}(X^{(n)},\mathbb{F}_5)$ into block-diagonal form, with one component for each irredicuble representation of $S_n$, noting that the coefficient field $\mathbb{F}_p$ satisfies $p>n$. 
The homology computation took very long for the $n=4$ case, several hours for each block in each dimension, though this could be made much more efficient 
by incorporating group actions into 
state of the art methods for 
homology calculations. We found for $n=4$ that
\[\sum_{k=0}^4 (-t)^k \ch H_k(X^{(4)},\mathbb{F}_5)=
\chi_{(4)}-(\chi_{4}+\chi_{3,1}+\chi_{2,2})t\]
\[+(2\chi_{3,1}+\chi_{2,2}+\chi_{2,1,1})t^2-(\chi_{3,1}+\chi_{2,1,1})t^3,\]
where $\chi_{\lambda}$ denotes the irreducible character of $S_4$ for a given Young diagram $\lambda$. This agrees with the predicted value for the full
space $\conf_n(\mathbb{R}^2)$ that one would obtain from equation \eqref{eq:lehrersolomon}. Replacing each irreducible character with their corresponding dimensions
\[(\chi_4,\chi_{3,1},\chi_{2,2},\chi_{2,1,1},\chi_{1,1,1,1})\mapsto (1,3,2,3,1),\]
we recover the expected value of
$1-6t+11t^2-6t^3$.

\bibliographystyle{plain}

\bibliography{refs}

\end{document}